\documentclass[oneside, english]{article}
\usepackage[T1]{fontenc}
\usepackage[latin1]{inputenc}
\usepackage{amssymb}
\usepackage{amsmath}
\usepackage[all]{xy}
\usepackage{amsthm}
\usepackage{setspace}

\makeatletter

  \theoremstyle{definition}
  \newtheorem{defn}{Definition}[section]

  \theoremstyle{remark}

   \theoremstyle{plain}
  \newtheorem{theorem}[defn]{Theorem}

    \newtheorem{example}[defn]{Example}
     
    \newcommand{\iso}{\overset{\sim}{\rightarrow}}

\usepackage{geometry}

\usepackage{babel}

\usepackage{amsmath,amssymb,amsfonts,latexsym}
\usepackage{mathrsfs}
\usepackage{epsf}
\usepackage{epsfig}
\usepackage[all]{xy}
\usepackage{enumerate}
\makeatother

\begin{document}

\def\AW[#1]^#2_#3{\ar@{-^>}@<.5ex>[#1]^{#2} \ar@{_<-}@<-.5ex>[#1]_{#3}}
 \def\NAW[#1]{\ar@{-^>}@<.5ex>[#1] \ar@{_<-}@<-.5ex>[#1]}

 \renewcommand{\arraystretch}{1}
 \renewcommand{\O}{\bigcirc}
 \newcommand{\OX}{\bigotimes}
 \newcommand{\OD}{\bigodot}
 \newcommand{\OV}{\O\llap{v\hspace{.6ex}}}
 \newcommand{\B}{\mbox{\Huge $\bullet$}}
  \newcommand{\D}{$\diamond$}

\xymatrixrowsep{1pc}
\xymatrixcolsep{1pc}

\pagestyle{myheadings}

\title{\textbf{On good $\mathbb{Z}$-gradings of basic Lie superalgebras}}
\author{Crystal Hoyt\footnote{Department of Mathematics, Bar-Ilan University, Ramat Gan 52900, Israel; hoyt@math.biu.ac.il.} \footnote{Supported by JSPS at Nara Women's University, Japan, and by ISF center of excellence 1691/10 at Bar-Ilan University, Israel. Supported by the Minerva foundation with funding from the Federal German Ministry for Education and Research.}}

\date{September 29, 2011}
\maketitle
\thispagestyle{empty}
\pagestyle{empty}

\vspace{-.5cm}
\begin{abstract}
We discuss the classification of good $\mathbb{Z}$-gradings of basic
Lie superalgebras. This problem arose in connection to $W$-algebras, where good $\mathbb{Z}$-gradings play a role in their construction.
\end{abstract}

\section{Introduction}
One component of the definition of a finite or affine super $W$-algebra is a good $\mathbb{Z}$-grading for a nilpotent element. Affine super $W$-algebras are (super) vertex algebras obtained from affine Lie superalgebras by quantum Hamiltonian reduction \cite{KRW}, whereas finite super $W$-algebras are associative superalgebras which can be defined via the universal enveloping algebra of a finite-dimensional simple Lie superalgebra \cite{DK}.

Good $\mathbb{Z}$-gradings of simple finite-dimensional Lie algebras where classified by A. Elashvili and V.G. Kac in 2005 \cite{EK}. K. Baur and N. Wallach classified nice parabolic subalgebras of reductive Lie algebras in \cite{BW}, which correspond to good even $\mathbb{Z}$-gradings by \cite[Theorem 2.1]{EK}.  J. Brundan and S. Goodwin classified good $\mathbb{R}$-gradings for semisimple finite-dimensional Lie algebras using certain polytopes, and proved that two finite $W$-algebras defined by the same nilpotent element $e\in\mathfrak{g}$ are isomorphic \cite{BG}.
This often allows one to reduce to the case that the good $\mathbb{Z}$-grading is even.

Here we discuss the classification of good $Z$-gradings of basic Lie superalgebras \cite{H11}.
In the case that $\mathfrak{g}$ is $\mathfrak{gl}(m|n)$ or $\mathfrak{osp}(m|2n)$ the good $\mathbb{Z}$-gradings are parameterized by ``good'' pyramids, generalizing the definition of \cite{EK}.  Whereas, for the exceptional Lie superalgebras, all good $\mathbb{Z}$-gradings are shown to be Dynkin.
Using this classification, one can determine which nilpotent elements have a good even $\mathbb{Z}$-grading.  For example, every nilpotent even element of $\mathfrak{gl}(m|n)$ has a good even $\mathbb{Z}$-grading.

\section{Basic Lie superalgebras}

Finite-dimensional simple Lie superalgebras were classified by V.G. Kac in \cite{K77}.  These can be separated into three types: basic, strange and Cartan.  A finite-dimensional simple Lie superalgebra  $\mathfrak{g} = \mathfrak{g}_{\bar{0}} \oplus\mathfrak{g}_{\bar{1}}$ is called {\em basic} if $\mathfrak{g}_{\bar{0}}$ is a reductive Lie algebra and $\mathfrak{g}$ has an even nondegenerate invariant bilinear form $(\cdot,\cdot)$.
This form is necessarily supersymmetric.  The basic Lie superalgebras are the following: $\mathfrak{sl}(m|n) : m\neq n$, $\mathfrak{psl}(n|n):=\mathfrak{sl}(n|n)/ \langle I_{2n} \rangle$,  $\mathfrak{osp}(m|2n)$, $D(2,1,\alpha)$, $F(4)$, $G(3)$, and finite dimensional simple Lie algebras.

Fix a Cartan subalgebra $\mathfrak{h}$. Then $\mathfrak{g}$ has a root space decomposition $\mathfrak{g}=\mathfrak{h}\oplus\bigoplus_{\alpha\in\Delta} \mathfrak{g}_{\alpha}.$  The $\mathbb{Z}/2\mathbb{Z}$-grading of $\mathfrak{g}$ determines a decomposition of  $\Delta$ into the disjoint union of the even roots $\Delta_{\bar{0}}$ and the odd roots $\Delta_{\bar{1}}$.  Corresponding to a set of simple roots
$\Pi=\{\alpha_1,\ldots,\alpha_n\}\subset\Delta$ of $\mathfrak{g}$, we have the triangular decomposition
$\mathfrak{g}=\mathfrak{n}^{-}\oplus\mathfrak{h}\oplus\mathfrak{n}^{+}$.

Most basic Lie superalgebras have more than one distinct Dynkin diagram.  This is due to the fact that the Weyl group does not act simply transitively on the set of bases.  However, we can extend the Weyl group to a Weyl groupoid by including ``odd reflections'', which allow us to move between the different bases.  In particular, if $\alpha_k\in\Pi$ is a simple isotropic root, then we can define the odd reflection at $\alpha_k$ to obtain a new set of simple roots $\Pi'$ for $\Delta$ \cite{LSS86}.

\section{Good $\mathbb{Z}$-gradings}

A $\mathbb{Z}$-grading $\mathfrak{g}=\oplus_{j\in\mathbb{Z}}\mathfrak{g}(j)$ is called {\em good} if there exists $e\in\mathfrak{g}_{\bar{0}}(2)$ such that the map $\mbox{ad }e:\mathfrak{g}(j)\rightarrow\mathfrak{g}(j+2)$ is injective for $j\leq -1$ and surjective for $j\geq -1$. If a $\mathbb{Z}$-grading of $\mathfrak{g}$ is defined by a semisimple element $h\in\mathfrak{g}_{\bar{0}}$, then this condition is equivalent to all of the eigenvalues of $\mbox{ad}(h)$ on the centralizer $\mathfrak{g}^{e}$ of $e$ in $\mathfrak{g}$ being non-negative.

An example of a good $\mathbb{Z}$-grading for a nilpotent element $e\in\mathfrak{g}_{\bar{0}}$ is the {\em Dynkin grading}. By the Jacobson-Morosov Theorem, $e$ belongs to an $\mathfrak{sl}_2$-triple $\mathfrak{s}=\{e,f,h\}\subset\mathfrak{g}_{\bar{0}}$, where $[e,f]=h$, $[h,e]=2e$ and $[h,f]=-2f$. By $\mathfrak{sl}_2$-theory, the grading of $\mathfrak{g}$ defined by $\mbox{ad }h$ is a good $\mathbb{Z}$-grading for $e$.

For each nilpotent even element $x\in\mathfrak{g}$ (up to conjugacy) we describe all $\mathbb{Z}$-gradings for which this element is good.  For the exceptional Lie superalgebras, we have the following

\begin{theorem}[Hoyt \cite{H11}]
All good $\mathbb{Z}$-gradings of the exceptional Lie superalgebras, $F(4)$, $G(3)$, and $D(2,1,\alpha)$, are Dynkin gradings.
\end{theorem}

To describe the good $\mathbb{Z}$-gradings of $\mathfrak{gl}(m|n)$ we generalize the definition of a pyramid given in \cite{BG, EK}.   A pyramid $P$ is a finite collection of boxes of size $2\times 2$ in the upper half plane which are centered at integer coordinates, such that for each $j=1,\ldots,N$, the second coordinates of the $j^{th}$ row equal $2j-1$ and the first coordinates of the $j^{th}$ row form an arithmetic progression  $f_j,f_j+2,\ldots,l_j$ with difference $2$, such that the first row is centered at $(0,0)$, i.e.  $f_1=-l_1$, and
\begin{equation}\label{eqpyramid} f_{j}\leq f_{j+1}\leq l_{j+1}\leq l_{j} \hspace{.5cm} \text{for all } j.\end{equation}
Each box of $P$ has even or odd parity.  We say that $P$ has {\em size} $(m|n)$ if $P$ has exactly $m$ even boxes and $n$ odd boxes.

Fix $m,n\in\mathbb{Z}_{+}$ and let $(p,q)$ be a partition of $(m|n)$. Let $r=\psi(p,q)\in Par(m+n)$ be the total ordering of the partitions $p$ and $q$ which satisfies: if $p_i=q_j$ for some $i,j$ then $\psi(p_i)< \psi(q_j)$. We define $Pyr(p,q)$ to be the set of pyramids which satisfy the following two conditions:
(1) the $j^{th}$ row of a pyramid $P\in Pyr(p,q)$ has length $r_j$;
(2) if $\psi^{-1}(r_j)\in p$ (resp.  $\psi^{-1}(r_j)\in q$) then all boxes in the $j^{th}$ row have even (resp. odd parity) and we mark these boxes with a ``$+$''  (resp. ``$-$'' sign).

Corresponding to each pyramid $P\in Pyr(p,q)$ we define a nilpotent element $e(P)\in\mathfrak{g}_{\bar{0}}$ and semisimple element $h(P)\in\mathfrak{g}_{\bar{0}}$, as follows.  Recall $\mathfrak{gl}(m|n)=\mbox{End}(V_{0}\oplus V_{1})$.  Fix a basis $\{v_1,\ldots,v_m\}$ of $V_{0}$ and $\{v_{m+1},\ldots,v_{m+n}\}$ of $V_{1}$.
Label the even (resp. odd) boxes of $P$ by the basis vectors of $V_{0}$ (resp. $V_{1}$). Define an endomorphism $e(P)$ of $V_{0}\oplus V_{1}$ as acting along the rows of the pyramid, i.e. by sending a basis vector $v_i$ to the basis vector which labels the box to the right of the box labeled by $v_i$ or to zero if it has no right neighbor.   Then $e(P)$ is nilpotent and corresponds to the partition $(p,q)$.  Since $e(P)$ does not depend the choice of $P$ in $Pyr(p,q)$, we may denote it by $e_{p,q}$.  Moreover, $e_{p,q}\in\mathfrak{g}_{\bar{0}}$ because boxes in the same row have the same parity.

Define $h(P)$ to be the $(m+n)$-diagonal matrix where the $i^{th}$ diagonal entry is the first coordinate of the box labeled by the basis vector $v_{i}$.  Then $h(P)$ defines a $\mathbb{Z}$-grading of $\mathfrak{g}$ for which $e_{p,q}\in\mathfrak{g}(2)$.  Let $P_{p,q}$ denote the symmetric pyramid from $Pyr(p,q)$. Then $h(P_{p,q})$ defines a Dynkin grading for $e_{p,q}$, and $P_{p,q}$ is called the {\em Dynkin pyramid} for the partition $(p|q)$.

\begin{theorem}[Hoyt \cite{H11}]  Let  $\mathfrak{g}=\mathfrak{gl}(m|n)$,
and let $(p,q)$ be a partition of $(m|n)$. If  $P$ is a pyramid from $Pyr(p,q)$, then the pair $(h(P),e_{p,q})$ is good. Moreover, every good grading for $e_{p,q}$ is of the form $(h(P),e_{p,q})$ for some pyramid $P\in Pyr(p,q)$.
\end{theorem}

This theorem is proven by studying the centralizer of a nilpotent element and of an $\mathfrak{sl}_2$ triple in $\mathfrak{gl}(m|n)$.  In a similar manner, we classify the good $\mathbb{Z}$-gradings for the Lie superalgebra $\mathfrak{osp}(m|2n)$ (see \cite{H11}).

A good $\mathbb{Z}$-grading of the Lie superalgebra $\mathfrak{g}$ for a nilpotent element $e\in\mathfrak{g}_{\bar{0}}$ restricts to a good $\mathbb{Z}$-grading for the Lie algebra $\mathfrak{g}_{\bar{0}}$.  So it is natural to ask which good $\mathbb{Z}$-gradings of $\mathfrak{g}_{\bar{0}}$ extend to a good $\mathbb{Z}$-grading of $\mathfrak{g}$, and to what extent is an extension unique.

\begin{example}\label{expyr1}
Let $\mathfrak{g}=\mathfrak{gl}(4|6)$ and consider the partitions $p=(3,1)$ and $q=(4,2)$.  The Dynkin grading of $\mathfrak{g}_{\bar{0}}=\mathfrak{gl}(4)\times\mathfrak{gl}(6)$ for the partition $(p,q)$ corresponds to the following symmetric pyramids.

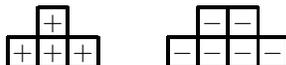
\setlength{\unitlength}{.4cm}
\begin{center}
\begin{picture}(5,2)(0,0)
\linethickness{1pt}
\put(0,0){\line(1,0){3}}
\put(0,1){\line(1,0){3}}
\put(1,2){\line(1,0){1}}
\put(0,0){\line(0,1){1}}
\put(1,0){\line(0,1){2}}
\put(2,0){\line(0,1){2}}
\put(3,0){\line(0,1){1}}
\put(0.15,0.25){$+$}
\put(1.15,0.25){$+$}
\put(2.15,0.25){$+$}
\put(1.15,1.25){$+$}
\end{picture}
\begin{picture}(6,2)(0,0)
\linethickness{1pt}
\put(0,0){\line(1,0){4}}
\put(0,1){\line(1,0){4}}
\put(1,2){\line(1,0){2}}
\put(0,0){\line(0,1){1}}
\put(1,0){\line(0,1){2}}
\put(2,0){\line(0,1){2}}
\put(3,0){\line(0,1){2}}
\put(4,0){\line(0,1){1}}
\put(0.15,0.25){$-$}
\put(1.15,0.25){$-$}
\put(2.15,0.25){$-$}
\put(3.15,0.25){$-$}
\put(1.15,1.25){$-$}
\put(2.15,1.25){$-$}
\end{picture}
\end{center}

There exist pyramids in $Pyr(p,q)$ for which the induced grading of $\mathfrak{g}_{\bar{0}}$ is the one given above, and these correspond to good $\mathbb{Z}$-gradings. They are represented by the following pyramids:\\
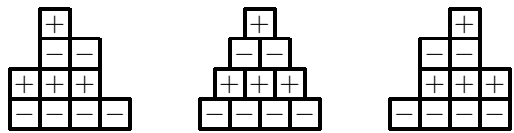
\begin{center}
\begin{picture}(6,4)(0,0)
\linethickness{1pt}
\put(0,0){\line(1,0){4}}
\put(0,1){\line(1,0){4}}
\put(0,2){\line(1,0){3}}
\put(1,3){\line(1,0){2}}
\put(1,4){\line(1,0){1}}
\put(0,0){\line(0,1){2}}
\put(1,0){\line(0,1){4}}
\put(2,0){\line(0,1){4}}
\put(3,0){\line(0,1){3}}
\put(4,0){\line(0,1){1}}
\put(0.15,1.25){$+$}
\put(1.15,1.25){$+$}
\put(1.15,3.25){$+$}
\put(2.15,1.25){$+$}
\put(0.15,0.25){$-$}
\put(1.15,0.25){$-$}
\put(2.15,0.25){$-$}
\put(3.15,0.25){$-$}
\put(1.15,2.25){$-$}
\put(2.15,2.25){$-$}
\end{picture}
\begin{picture}(6,4)(0,0)
\linethickness{1pt}
\put(0,0){\line(1,0){4}}
\put(0,1){\line(1,0){4}}
\put(.5,2){\line(1,0){3}}
\put(1,3){\line(1,0){2}}
\put(1.5,4){\line(1,0){1}}
\put(0,0){\line(0,1){1}}
\put(1,0){\line(0,1){1}}
\put(2,0){\line(0,1){1}}
\put(3,0){\line(0,1){1}}
\put(4,0){\line(0,1){1}}
\put(.5,1){\line(0,1){1}}
\put(1.5,1){\line(0,1){1}}
\put(2.5,1){\line(0,1){1}}
\put(3.5,1){\line(0,1){1}}
\put(1,2){\line(0,1){1}}
\put(2,2){\line(0,1){1}}
\put(3,2){\line(0,1){1}}
\put(1.5,3){\line(0,1){1}}
\put(2.5,3){\line(0,1){1}}
\put(0.65,1.25){$+$}
\put(1.65,1.25){$+$}
\put(2.65,1.25){$+$}
\put(1.65,3.25){$+$}
\put(0.15,0.25){$-$}
\put(1.15,0.25){$-$}
\put(2.15,0.25){$-$}
\put(3.15,0.25){$-$}
\put(1.15,2.25){$-$}
\put(2.15,2.25){$-$}
\end{picture}
\begin{picture}(6,4)(0,0)
\linethickness{1pt}
\put(0,0){\line(1,0){4}}
\put(0,1){\line(1,0){4}}
\put(1,2){\line(1,0){3}}
\put(1,3){\line(1,0){2}}
\put(2,4){\line(1,0){1}}
\put(0,0){\line(0,1){1}}
\put(1,0){\line(0,1){3}}
\put(2,0){\line(0,1){4}}
\put(3,0){\line(0,1){4}}
\put(4,0){\line(0,1){2}}
\put(1.15,1.25){$+$}
\put(2.15,1.25){$+$}
\put(3.15,1.25){$+$}
\put(2.15,3.25){$+$}
\put(0.15,0.25){$-$}
\put(1.15,0.25){$-$}
\put(2.15,0.25){$-$}
\put(3.15,0.25){$-$}
\put(1.15,2.25){$-$}
\put(2.15,2.25){$-$}
\end{picture}
\end{center}
\end{example}

\begin{example}\label{expyr2}
Let $\mathfrak{g}=\mathfrak{gl}(4|6)$ and consider the partitions $p=(3,1)$ and $q=(4,2)$.
The following pyramids represent a good $\mathbb{Z}$-grading of $\mathfrak{g}_{\bar{0}}$ for which
there is no good $\mathbb{Z}$-grading of $\mathfrak{g}$ with this induced good $\mathbb{Z}$-grading of $\mathfrak{g}_{\bar{0}}$.\\
\setlength{\unitlength}{.4cm}
\begin{center}
\begin{picture}(5,2)(0,0)
\linethickness{1pt}
\put(0,0){\line(1,0){3}}
\put(0,1){\line(1,0){3}}
\put(0,2){\line(1,0){1}}
\put(0,0){\line(0,1){2}}
\put(1,0){\line(0,1){2}}
\put(2,0){\line(0,1){1}}
\put(3,0){\line(0,1){1}}
\put(0.15,0.25){$+$}
\put(1.15,0.25){$+$}
\put(2.15,0.25){$+$}
\put(0.15,1.25){$+$}
\end{picture}
\begin{picture}(6,2)(0,0)
\linethickness{1pt}
\put(0,0){\line(1,0){4}}
\put(0,1){\line(1,0){4}}
\put(2,2){\line(1,0){2}}
\put(0,0){\line(0,1){1}}
\put(1,0){\line(0,1){1}}
\put(2,0){\line(0,1){2}}
\put(3,0){\line(0,1){2}}
\put(4,0){\line(0,1){2}}
\put(0.15,0.25){$-$}
\put(1.15,0.25){$-$}
\put(2.15,0.25){$-$}
\put(3.15,0.25){$-$}
\put(2.15,1.25){$-$}
\put(3.15,1.25){$-$}
\end{picture}
\end{center}
\end{example}

\section{Centralizers of $\mathfrak{sl}_{2}$-triples}

The centralizers of $\mathfrak{sl}_{2}$-triples in $\mathfrak{gl}(m|n)$ and $\mathfrak{osp}(m|2n)$ can be described following the ideas of \cite{J} for the Lie algebras $\mathfrak{gl}(m)$, $\mathfrak{so}(m)$ and $\mathfrak{sp}(2n)$.
There is a one-to-one correspondence between $G$-orbits of nilpotent even elements in $\mathfrak{gl}(m|n)$ and partitions of $(m|n)$.  Let $p=(r_1^{m_{1}},\ldots,r_{N}^{m_N})$ be a partition of $m$ and $q=(r_1^{n_{1}},\ldots,r_{N}^{n_N})$ a partition of $n$, that is $r_i$ has multiplicity $m_i$ in $p$ and multiplicity $n_i$ in $q$.  We note that $m_i$ or $n_i$ may be zero.

\begin{theorem}[Hoyt \cite{H11}]\label{thmgl} Let $\mathfrak{g}=\mathfrak{gl}(m|n)$. Let $e$ be a nilpotent even element corresponding to a partition $(p,q)$ of $(m|n)$,and let $\mathfrak{s}=\{e,f,h\}\subset\mathfrak{g}_{\bar{0}}'$ be an $\mathfrak{sl}_2$-triple for $e$. Then we have an isomorphism
$$\mathfrak{g}^{\mathfrak{s}}\iso \mathfrak{gl}(m_1,n_1)\times  \cdots\times \mathfrak{gl}(m_N,n_N)$$ of Lie superalgebras.
\end{theorem}

A partition is called {\em  symplectic} (resp. {\em orthogonal}) if $m_{p_i}$ is even for odd $p_i$ (resp. even $p_i$).
We say that a partition $(p,q)$ of $(m|2n)$ is orthosymplectic if $p$ is an orthogonal partition of $m$ and $q$ is a symplectic partition of $2n$. There is a one-to-one correspondence between $G$-orbits of nilpotent even elements in $\mathfrak{osp}(m|2n)$ and orthosymplectic partitions of $(m|2n)$.  Let $(p,q)$ be an orthosymplectic partition of $(m|2n)$, and represent it as $p=(r_1^{m_{1}},\ldots,r_{N}^{m_N},s_1^{2c_{1}},\ldots,s_{T}^{2c_{T}})$ and $q=(r_1^{2n_{1}},\ldots,r_{N}^{2n_N},s_1^{d_{1}},\ldots,s_{T}^{d_T})$, where $r_i$ are the even parts and $s_i$ are the odd parts.

\begin{theorem}[Hoyt \cite{H11}] Let $\mathfrak{g}=\mathfrak{osp}(m|2n)$.
Let $e$ be a nilpotent even element corresponding to an orthosymplectic partition $(p,q)$ of $(m|n)$, and let $\mathfrak{s}=\{e,f,h\}\subset\mathfrak{g}_{\bar{0}}'$ be an $\mathfrak{sl}_2$-triple for $e$. Then we have an isomorphism
$$\mathfrak{g}^{\mathfrak{s}}\iso \mathfrak{osp}(m_1,2n_1)\times \cdots  \times \mathfrak{osp}(m_N,2n_N)\times\mathfrak{osp}(d_1,2c_1)\times   \cdots  \times \mathfrak{osp}(d_T,2c_T)$$
of Lie superalgebras.
\end{theorem}

\section{$\mathbb{Z}$-gradings and the Weyl groupoid}

Let $\mathfrak{g}$ be a basic Lie superalgebra, $\mathfrak{g}\neq\mathfrak{psl}(2|2)$, or let $\mathfrak{g}$ be $\mathfrak{gl}(m|n)$ or $\mathfrak{sl}(n|n)$. Fix a $\mathbb{Z}$-grading  $\mathfrak{g}=\oplus_{j\in\mathbb{Z}} \mathfrak{g}(j)$ satisfying $\mathfrak{Z}(\mathfrak{g}_{\bar{0}})\subset \mathfrak{g}_{\bar{0}}(0)$.
Fix a Cartan subalgebra $\mathfrak{h}\subset\mathfrak{g}_{\bar{0}}(0)$, and let $\Delta$ be the set of roots.
The root space decomposition with respect to $\mathfrak{h}$ is compatible with the $\mathbb{Z}$-grading, so we can define a map $\mbox{Deg}:\Delta\cup\{0\}\rightarrow\mathbb{Z}$ by $\mbox{Deg}(\alpha)= k$ if $\alpha\in\Delta_k$ and $\mbox{Deg}(0)=0$.

Now for each base $\Pi\subset\Delta$, the degree map of a $\mathbb{Z}$-grading is determined by its restriction to $\Pi$, that is, by $D:\Pi\rightarrow\mathbb{Z}$. A reflection at a simple root of $\Pi$ yields a new map $D':\Pi'\rightarrow\mathbb{Z}$, where $\Pi'$ is the reflected base and $D'$ is defined on $\Pi'$ by linearity.  The maps $D:\Pi\rightarrow\mathbb{Z}$ and  $D':\Pi'\rightarrow\mathbb{Z}$ define the same grading.

It is natural to ask the following question:
when do two maps $D_1:\Pi_1\rightarrow\mathbb{N}$ and $D_2:\Pi_2\rightarrow\mathbb{N}$ define the same $\mathbb{Z}$-grading, i.e. when can they be extended to a linear map $\mbox{Deg}:\Delta\cup\{0\}\rightarrow\mathbb{Z}$?

\begin{theorem}[Hoyt \cite{H11}]
 Let $\Pi_1=\{\alpha_1,\ldots,\alpha_n\},\ \Pi_2=\{\beta_1,\ldots,\beta_n\}$ be two different bases for $\Delta$.
The maps $D_1:\Pi_1\rightarrow\mathbb{N}$, $D_2:\Pi_2\rightarrow\mathbb{N}$ define the same grading if and only if there is a sequence of even and odd reflections $\mathcal{R}$ at simple roots of degree zero such that (after reordering) $\mathcal{R}(\alpha_i)=\beta_i$ and $D_1(\alpha_i)=D_2(\beta_i)$ for $i=1,\ldots,n$.

Two Dynkin diagrams $\Gamma_1,\Gamma_2$ for a basic Lie superalgebra $\mathfrak{g}$ with degree maps $D_i:\Gamma_i\rightarrow\mathbb{N}$ define the same $\mathbb{Z}$-grading if and only if there is a sequence of odd reflections $\mathcal{R}$ at simple isotropic roots of degree zero such that $\mathcal{R}(\Gamma_1)=\Gamma_2$ and $D_1=D_2$ with the ordering of the vertices defined by $\mathcal{R}$.  This defines an equivalence relation on Dynkin diagrams with nonnegative integer labels.

\end{theorem}

\end{document}